\documentclass{elsart}

\usepackage{amssymb}

\usepackage[latin1]{inputenc}
\usepackage[centertags]{amsmath}
\usepackage{amsfonts}
\usepackage{amsmath}
\usepackage{amssymb}
\usepackage{latexsym}

\usepackage{mathrsfs}
\usepackage{newlfont}

\newcommand\CC{{\mathbb C}}

\newcommand{\D}{\displaystyle}

\begin{document}

\begin{frontmatter}



\title{Some examples of orthogonal matrix polynomials satisfying odd order differential equations}


\author{Antonio J. Durán and Manuel D. de la Iglesia}

\address{Departamento de An\'{a}lisis Matem\'{a}tico.
   Universidad de Sevilla \\
   \footnotesize Apdo (P. O. BOX) 1160. 41080 Sevilla. Spain.\\
   duran@us.es, mdi29@us.es}

\begin{abstract}
It is well known that if a finite order linear differential operator
with polynomial coefficients has as eigenfunctions a sequence of
orthogonal polynomials with respect to a positive measure (with
support in the real line), then its order has to be even. This
property no longer holds in the case of orthogonal matrix
polynomials. The aim of this paper is to present examples of
weight matrices such that the corresponding sequences of  matrix
orthogonal polynomials are eigenfunctions of certain linear
differential operators of odd order. The weight matrices are of the
form
$$
W(t)=t^{\alpha}e^{-t}e^{At}t^{B}t^{B^*}e^{A^* t},
$$
where $A$ and $B$ are certain (nilpotent and diagonal, respectively)
$N\times N$ matrices. These weight matrices are the first examples
illustrating this new phenomenon which are not reducible to scalar
weights.
\end{abstract}

\end{frontmatter}

\section{Introduction}\label{int}

During the last few years many families of orthogonal matrix
polynomials have been found that are eigenfunctions of some fixed
second order linear differential operator with matrix coefficients
which do not depend on the degree of the polynomial. When the
corresponding eigenvalues are Hermitian then the differential
operator for the orthonormal polynomials $(P_n)_n$ is also symmetric
with respect to (the inner product defined from) its orthogonalizing
weight matrix $W$, and conversely. More precisely, for an operator
(whose coefficients are multiplied on the right)
\begin{equation}\label{difope}
\ell_2=D^2A_2(t)+D^1A_1(t)+D^0A_0,
\end{equation}
and a family of orthonormal polynomials $(P_n)_n$ with respect to a
weight matrix $W$, the following conditions are equivalent:
\begin{enumerate}
\item $\ell _2(P_n(t))=\Gamma _nP_n(t)$,
with $\Gamma_n$ Hermitian
\item $\int \ell _2 (P)dWQ^*=\int PdW\ell _2(Q)^*$, \quad $P,Q\in \CC^{N\times N}[t]$.
\end{enumerate}

An overview of those families of orthogonal matrix polynomials
satisfying second order differential operators shows some important
differences with respect to the scalar situation. First of all, it
has been necessary to develop new techniques (rather different from
the scalar ones) to find examples (see \cite{DG1}, \cite{D2},
\cite{GPT1}, \cite{GPT2} and \cite{GPT3}). Secondly, the complexity
of the matrix world has led to an embarrassment of riches if we
compare with the only scalar examples of Hermite, Laguerre and
Jacobi. Also some new phenomena have appeared. For instance, some
families of orthogonal polynomials $(P_n)_n$ have been found which
are common eigenfunctions of several linearly independent second
order differential operators (see \cite{CG2}, \cite{dIG} or
\cite{DL}). This shows that the algebras of differential operators
having these families of orthogonal matrix polynomials as
eigenfunctions are going to have a richer structure than the
corresponding algebras in the scalar case, where for the classical
scalar families those algebras reduce to the associated second order
differential operator and any polynomial in this operator (see
\cite{M1}).

The purpose of this paper is to show a new phenomenon. Indeed, in
the scalar case it is well known that if a differential operator has
a sequence of orthogonal polynomials as eigenfunctions then its
order has to be even (see \cite{K}). This situation is not true in
the matrix case. We present here some examples of orthogonal matrix
polynomials which are eigenfunctions of certain differential
operators of odd order. Those examples are particular cases of
orthogonal matrix polynomials with respect to a wider family of
$N\times N$ weight matrices defined by
\begin{equation}\label{W}
W_{\alpha
,\nu_1,\cdots,\nu_{N-1}}(t)=t^{\alpha}e^{-t}e^{At}t^{\frac{1}{2}J}t^{\frac{1}{2}J^*}e^{A^*
t}, \quad \alpha>-1,\quad t\in (0,+\infty),
\end{equation}
where $A$ is the $N\times N$ nilpotent matrix:
\begin{equation*}
A=
\begin{pmatrix}
    0 & \nu_1 & 0  & \cdots & 0 \\
    0 & 0 & \nu_2  & \cdots & 0 \\
    \vdots & \vdots & \vdots  & \ddots & \vdots \\
    0 & 0 & 0  & \cdots & \nu_{N-1} \\
    0 & 0 & 0  & \cdots & 0 \\
  \end{pmatrix} ,\quad \nu_i\in\mathbb{C}\setminus\{0\},\quad i=1,\cdots,N-1,
\end{equation*}
and $J$ is the $N\times N$ diagonal matrix:
\begin{equation*}
J=\begin{pmatrix}
    N-1 & 0   & \cdots& 0 & 0 \\
    0 & N-2   & \cdots& 0 & 0 \\
    \vdots & \vdots   & \ddots& \vdots & \vdots \\
    0 & 0   & \cdots & 1& 0 \\
    0 & 0   & \cdots& 0 & 0 \\
  \end{pmatrix} .
\end{equation*}
The weight matrix corresponding to $N=2$ yields the following third
order symmetric differential operator:
\begin{align}\label{ope3}
\nonumber \ell _{3,1}=&D^3\begin{pmatrix}
                     -|a|^2t^2 & at^2(1+|a|^2t) \\
                     -\bar{a}t & |a|^2t^2 \\
                   \end{pmatrix}
    +\\ \nonumber &D^2 \begin{pmatrix}
    -t(2+|a|^2(\alpha+5)) & at(2\alpha+4+t(1+|a|^2(\alpha+5))) \\
    -\bar{a}(\alpha+2) & t(2+|a|^2(\alpha+2)) \\
    \end{pmatrix}+\\
    \nonumber &D^1\begin{pmatrix}
    t-2(\alpha+2)(1+|a|^2) & \begin{Large}\mbox{$\frac{|a|^2(\alpha+1)(\alpha+2)+t(1+2|a|^2(1+|a|^2(\alpha+2))}{\bar{a}}$}\end{Large} \\
                                                       -\D\frac{1}{a} & 2\alpha+2-t \\
                                                          \end{pmatrix}+\\
                                         &D^0\begin{pmatrix}
                                         1+\alpha & -\D\frac{1}{\bar{a}}(1+\alpha)(|a|^2\alpha-1) \\
                                         \D\frac{1}{a} & -(1+\alpha) \\
                                         \end{pmatrix}.
\end{align}

In fact, the weight matrix has two linearly independent symmetric
third order differential operators and other two linearly independent operators of order two.
The weight matrices corresponding to $N=3$ and $4$ yield to fifth and
seventh order (respectively) symmetric differential operator, but no
third order. This is related to the fact that $A$ is nilpotent of
order $N$. Our conjecture is that each $W_{\alpha
,\nu_1,\cdots,\nu_{N-1}}$ has a symmetric differential operator of
order $2N-1$, as well as other of even order (as we will see below,
$W_{\alpha ,\nu_1,\cdots,\nu_{N-1}}$ always has  at least one
symmetric second order differential operator).

Some other examples of $2\times2$ situations have recently appeared
in the literature having differential operators of odd order (order
1, to be more precise): see \cite{CG1}, \cite{CG2} and \cite{DG2}.
Some of them correspond to positive definite weight matrices and
happen to be symmetric. Some others do not even correspond to a
positive definite weight matrix. However, the first ones have
weights that reduce to scalar weights. We say that a weight matrix
$W$ reduces to scalar weights if there exists a nonsingular matrix
$T$ independent of $t$ for which
$$
W(t)=TD(t) T^*,
$$
with $D(t)$ diagonal. Weight matrices reducible to scalar weights
are actually a collection of $N$ scalar weights. When they are
considered by themselves and not in connection with differential
equations these weights belong to the study of scalar orthogonality
more than to the matrix one. The $2\times 2$ positive definite
weight matrices referred above (see \cite{CG1}, \cite{CG2} and
\cite{DG2}) are of the form $W(t)=T\begin{pmatrix} \omega_1(t)&0 \\
0& \omega_2(t)\end{pmatrix}T^*$, where $T$ is a nonsingular matrix
(independent of $t$). We have to point out that, however, this
matrix $T$ do not factorize their associated differential operators
(that is, that operators are not of the form $T\begin{pmatrix} d_1&0
\\ 0& d_2\end{pmatrix}T^*$; see \cite{GPT2} for a notion of scalar
reducibility for the pair consisting of the weight and the
differential operator). According to the previous definition, the
examples shown in this paper do not reduce to scalar weight and
hence are the first weight matrices with this property that have a
symmetric differential operator of odd order.

The paper is organized as follows. In Section \ref{sect2}, we
introduce the family $W_{\alpha ,\nu_1,\cdots,\nu_{N-1}}$ (see
(\ref{W})) and study its symmetric second order differential
operators. We find that
$$
\ell _{2,1}=D^2tI+D^1[(\alpha +1)I+J+t(A-I)]+D^0 [(J+\alpha I)A-J]
$$
is symmetric for $W_{\alpha ,\nu_1,\cdots,\nu_{N-1}}$. When the
moduli of the parameters $\nu_i$, $i=1,\cdots, N-2$, are defined
from $\nu_{N-1}$ by
$$
i(N-i)|\nu_{N-1}|^2=(N-1)|\nu_{i}|^2+(N-i-1)|\nu_{i}|^2|\nu_{N-1}|^2,
$$
our weight matrix enjoys another symmetric second order differential
operator $\ell _{2,2}$ with leading coefficient $A_2(t)=t(J-At)$. It
turns out that the operators $\ell _{2,1}$ and $\ell _{2,2}$ commute
and satisfy the polynomial equation
$$
\prod_{i=1}^N
\bigg((i-1)\ell_{2,1}-\ell_{2,2}+\left[\frac{(N-1)(N-i)}{|\nu_{N-1}|^2}+(i-1)(N-i)\right]I\bigg)=0.
$$
We note here that the symbol 0 is used to denote either zero scalar
or zero matrix (all entries equal to zero scalar) while $I$ is used
to denote the identity matrix of dimension determined from the
context.

We complete Section \ref{sect2} showing some structural formulas
(Rodrigues' formula and three--term recurrence relation) for a
sequence $(\mathcal{P}_{n,\alpha,a})_n$ of orthogonal polynomials
with respect to the weight matrix
\begin{align}\label{ejemplo}
W_{\alpha, a}(t)=t^\alpha e^{-t}\begin{pmatrix} t(1+\vert a\vert
^2t) &at\\ \bar a t& 1\end{pmatrix}, \quad \alpha >-1, t>0,
\end{align}
which is the case $N=2$ of (\ref{W}) (putting $\nu_1=a$). This
weight matrix appears for the first time in \cite{CMV}, where the
authors found a Pearson equation for it and proved that the
derivatives $(\mathcal{P'}_{n,\alpha, a})_n$ are again orthogonal
with respect to a certain weight matrix.

In order to look for symmetric differential operators of higher
order with respect to the weight matrix $W_{\alpha
,\nu_1,\cdots,\nu_{N-1}}$, we show, in Section \ref{Hodo}, how the
symmetry of a differential operator $\ell _k$ of any order $k$ with
respect to a weight matrix $W$ can be reduced to a set of $k+1$
differential equations (of order $0, 1, \cdots , k$, respectively)
and certain boundary conditions for $W$ and the coefficients of
$\ell _k$ at the endpoints of the support of $W$. This can be seen
as the matrix valued version of Littlejohn's conditions for
the symmetry over polynomials of this kind of differential operators
in the scalar case (see \cite{L}).

Finally, in Section \ref{alg}, we study the algebra of differential
operators associated with (\ref{ejemplo}) defined by
$$
\mathcal{D}(W_{\alpha ,a})=\left\{\ell =\sum_{i=0}^k D^i A_i(t) :
\ell (\mathcal
P_{n,\alpha,a}(t))=\Gamma_n(\ell)\mathcal{P}_{n,\alpha
,a}(t),\;n\geq0\right\}.
$$
In particular, we find two linearly independent symmetric third
order differential operators (one of them is given in (\ref{ope3})).
It seems that two new linearly independent operators appear as one
increases by one the order of the operators in question. We give
strong computational evidences which support our conjecture that the
algebra $\mathcal{D}(W_{\alpha,a})$ is generated by the set $\{ I,
\ell _{2,1}, \ell_{3,1} \}$ (except for some exceptional values
involving the parameters $\alpha$ and $a$).

\section{Second order differential operators for $W_{\alpha
,\nu_1,\cdots,\nu_{N-1}}$}\label{sect2}

The aim of this section is to introduce the weight matrix $W_{\alpha
,\nu_1,\cdots,\nu_{N-1}}$ and study its second order differential
operators.

Under the assumption that
\begin{equation*}
A_2(t)W(t) \mbox{ and } (A_2(t)W(t))'- A_1(t)W(t),
\end{equation*}
vanish at each of the endpoints of the support of $W(t)$, to find
symmetric second order differential operators for a weight matrix
$W$ is enough to solve the  following equations:
\begin{align}\label{eq2.1}
A_2W&=WA_2^*,\\\label{eq2.2} 2(A_2W)'&=WA_1^*+A_1W, \\\label{eq2.3}
(A_2W)''-(A_1W)'+A_0W&=WA_0^*,
\end{align}
(see Theorem 3.1 of \cite{DG1} or also \cite{GPT2}).

Assuming that $A_2(t)$ is a scalar matrix,  it has been proved in
\cite{DG1} that the differential equation (\ref{eq2.2}) is
equivalent to the fact that $W$ can be factorized in the form
$W(t)=\rho (t)T(t)T^*(t)$, where $\rho $ is a scalar function and
$T$ is a matrix function satisfying certain first order differential
equation. When $\rho(t) =t^\alpha e^{-t}$, $\alpha>-1$,  i.e. the
Laguerre classical scalar weight, and $A_2(t)=tI$, this first order
differential equation for $T$ takes the form
\begin{equation}\label{Tp}
T'(t)=\left( A+\frac{B}t\right)T(t).
\end{equation}
Let us consider a weight matrix of the form
$W(t)=t^{\alpha}e^{-t}T(t)T(t)^*$, where
\begin{equation}\label{T}
T(t)=e^{At}t^{B}=e^{At}e^{B \log t},
\end{equation}
with $A$ and $B$ any matrices. Using the formula
\begin{equation*}
e^{At}B=\bigg(\sum_{n\geq0}\frac{t^n}{n!}\mbox{ad}_A^n B
\bigg)e^{At},
\end{equation*}
we can write the derivative of (\ref{T}) as
\begin{equation*}
T'(t)=AT+\frac{1}{t}TB=\bigg(\frac{B}{t}+A+\mbox{ad}_A
B+\sum_{n\geq2}\frac{t^{n-1}}{n!}\mbox{ad}_A^n B\bigg)T.
\end{equation*}
We use the standard notation
$$
\mbox{ad}_X^0 Y=Y, \quad \mbox{ad}_X^1 Y=[X,Y], \quad \mbox{ad}_X^2
Y=[X,[X,Y]],
$$
and, in general, $\mbox{ad}_X^{n+1} Y=[X,\mbox{ad}_X^n Y]$, where
$[X,Y]=XY-YX$.

In order for $T$ to satisfy a differential equation like (\ref{Tp}),
we need to choose matrices $A$ and $B$ for which $\mbox{ad}_A^2B=0$.
Once we have chosen the matrix
\begin{equation}\label{AA}
A= \begin{pmatrix}
    0 & \nu_1 & 0  & \cdots & 0 \\
    0 & 0 & \nu_2  & \cdots & 0 \\
    \vdots & \vdots & \vdots  & \ddots & \vdots \\
    0 & 0 & 0  & \cdots & \nu_{N-1} \\
    0 & 0 & 0  & \cdots & 0 \\
  \end{pmatrix} ,\quad \nu_i\in\mathbb{C}\setminus\{0\},\quad i=1,\cdots,N-1,
\end{equation}
and $B$ to be diagonal and singular, it follows from the condition
$\mbox{ad}_A^2B=0$ that necessarily $B=uJ$ where
\begin{equation}\label{JJ}
J=\begin{pmatrix}
    N-1 & 0   & \cdots& 0 & 0 \\
    0 & N-2   & \cdots& 0 & 0 \\
    \vdots & \vdots   & \ddots& \vdots & \vdots \\
    0 & 0   & \cdots & 1& 0 \\
    0 & 0   & \cdots& 0 & 0 \\
  \end{pmatrix} ,
  \end{equation}
and $u$ is any complex number. However, it turns out that for the
existence of a symmetric second order differential operator with
$A_2(t)=tI$, $u$ has to be $1/2$, except for $N=2$ where $u$
can be any complex number.

This is the reason why we have chosen the weight matrix to be
$$
W_{\alpha
,\nu_1,\cdots,\nu_{N-1}}(t)=t^{\alpha}e^{-t}e^{At}t^{\frac{1}{2}J}t^{\frac{1}{2}J^*}e^{A^*
t}, \quad \alpha>-1,\quad t\in (0,+\infty),
$$
where $A$ and $J$ are defined by (\ref{AA}) and (\ref{JJ}),
respectively, as a candidate to have symmetric second order
differential operators like (\ref{difope}) with $A_2(t)=tI$. Our
choice allows to write $W_{\alpha
,\nu_1,\cdots,\nu_{N-1}}(t)=t^{\alpha}e^{-t}T(t)T(t)^*$, where
$T(t)=e^{At}t^{\frac{1}{2}J}$ satisfies
\begin{equation}\label{Tp2}
T'(t)=\frac{1}{2}\bigg(\frac{J}{t}+A\bigg)T(t).
\end{equation}
Actually $W_{\alpha ,\nu_1,\cdots,\nu_{N-1}}$ is not a bad
candidate!

\begin{thm}
The second order differential operator
\begin{align}\label{difopeW}
\ell _{2,1}=D^2tI+D^1[(\alpha+1)I+J+t(A-I)]+D^0[(J+\alpha I)A-J],
\end{align}
is symmetric with respect to $W_{\alpha ,\nu_1,\cdots,\nu_{N-1}}$.
\end{thm}
\pf We only need to verify equation (\ref{eq2.3}). Using
(\ref{eq2.2}), that equation is equivalent to
$(A_1W-WA_1^*)'=2(A_0W-WA_0^*)$ (where, to simplify the notation, we
remove the dependence on $\alpha$ and $\nu_1,\cdots, \nu_{N-1}$).
Using the formulas
$$
e^{At}J=Je^{At}-Ate^{At}\quad\mbox{and}\quad
e^{A^*t}J=Je^{A^*t}+A^*te^{A^*t}
$$
we obtain
\begin{align*}
JW-WJ=&t(AW-WA^*)\\
JAW-WA^*J=&(AW-WA^*)+(AJW-WJA^*).
\end{align*}
Also, as a consequence of the first equation above, we deduce
\begin{align*}
JW'-W'J&=AW-WA^*+t(AW'-W'A^*)\\
t(A^2W-W(A^*)^2)&=(JWA^*-AWJ)+(AJW-WJA^*).
\end{align*}
Hence, the equation $(A_1W-WA_1^*)'=2(A_0W-WA_0^*)$ follows from these
equalities taking into account the differential equation for $T(t)$
in (\ref{Tp2}). \hfill $\Box$

For $N=2$, the weight matrix
$W(t)=t^{\alpha}e^{-t}e^{At}t^{B}t^{B^*}e^{A^*t}$, where
$$
A=\begin{pmatrix} 0 & v\\ 0 & 0\end{pmatrix},\quad B=\begin{pmatrix}
u&0\\ 0&0\end{pmatrix},
$$
has associated the symmetric second order differential operator
$$
\ell _2=D^2\begin{pmatrix} t&0\\ 0&t\end{pmatrix}+D^1\begin{pmatrix}
2u+\alpha+1-t&2tv(1-u)\\ 0&\alpha+1-t\end{pmatrix}
+D^0\begin{pmatrix} -1&v(1+\alpha)\\ 0&0\end{pmatrix}.
$$
As we wrote above for $N=3$, in order to have a symmetric second
order differential operator with $A_2(t)=tI$ it is necessary that
$B=\frac{1}{2}J$.

We now search for another second order differential operator $$\ell
_{2,2}=D^2A_2(t)+D^2A_1(t)+D^0A_0$$ for the weight matrix $W_{\alpha
,\nu_1,\cdots,\nu_{N-1}}$. We follow the lines of the method
developed in \cite{D2}. This method consists of looking for {\sl
good factorizations} of $W_{\alpha ,\nu_1,\cdots,\nu_{N-1}}$ in the
form $t^\alpha e^{-t} R(t)R^*(t)$. Under the assumption that
$A_2W_{\alpha ,\nu_1,\cdots,\nu_{N-1}}$ is Hermitian, by a good
factorization for $W_{\alpha ,\nu_1,\cdots,\nu_{N-1}}$ we mean that
the factor $R$ satisfies the first order differential equation
$R'=FR$, where
\begin{enumerate}
  \item The differential coefficient $F$ is related to $A_2$ and $A_1$ by the equation
\begin{align}\label{eqnm1.1}
A_1=A_2F+FA_2+C,
\end{align}
where $C(t)=(t^\alpha e^{-t} A_2(t))'/t^\alpha e^{-t}$ and
  \item The matrix function
\begin{align}\label{eqnm1.2}
R^{-1}(FA_2F+F'A_2+FC-A_0)R
\end{align}
is Hermitian.
\end{enumerate}
Then (1) guarantees that $W_{\alpha ,\nu_1,\cdots,\nu_{N-1}}$
satisfies equation (\ref{eq2.2}) and (2) equation (\ref{eq2.3}) (see
\cite{D2}, Section 2).

According to this approach, first of all, we have to look for the
leading coefficient $A_2$ of the differential operator $\ell
_{2,2}$. This coefficient has to satisfy $A_2W=WA_2^*$ (we are again
removing the dependence on $\alpha$ and $\nu_1,\cdots, \nu_N$). The
relation $AJ-JA=-A$ gives a very natural candidate for $A_2$.
Indeed, $A_2W=WA_2^*$ is equivalent to
$t^{-\frac{1}{2}J}e^{-At}A_2e^{At}t^{\frac{1}{2}J}$ being Hermitian.
If we put
\begin{equation}\label{A2}
    A_2(t)=t(J-At),
\end{equation}
a straightforward computation gives
\begin{align*}
   t^{-\frac{1}{2}J}e^{-At}(t(J-At))e^{At}t^{\frac{1}{2}J}=&t^{-\frac{1}{2}J}\bigg(\sum_{n\geq0}
   \frac{(-1)^nt^{n+1}}{n!}\mbox{ad}_A^n
   J-At^2\bigg)t^{\frac{1}{2}J}\\
   =&t\cdot t^{-\frac{1}{2}J}Jt^{\frac{1}{2}J}=tJ,
\end{align*}
implying that (\ref{eq2.1}) is satisfied since $tJ$ is Hermitian.

Once we have a candidate for $A_2$, we need to choose a certain good
factorization of the  weight matrix $W(t)=t^\alpha e^{-t}
R(t)R^*(t)$ (in the sense explained above). We proceed by taking a
unitary matrix function $U(t)$, and writing
$$
R(t)=e^{At}t^{\frac{1}{2}J}U(t),
$$
so that $W(t)=t^{\alpha}e^{-t}R(t)R(t)^*$. The matrix function
$U(t)$ will be introduced later.

The definition of $R(t)$ implies that the coefficient $F$ of the
first order differential equation $R'(t)=F(t)R(t)$ for $R$ is
\begin{equation}\label{F}
F(t)=\frac{1}{2}\bigg(\frac{1}{t}J+A\bigg)+e^{At}t^{\frac{1}{2}J}X(t)t^{-\frac{1}{2}J}e^{-At},
\end{equation}
where we have written $X(t)=U'(t)U(t)^{-1}$. Since $U(t)$ is
unitary, it turns out that the matrix function $X(t)$ is
skew-Hermitian. Taking this into account we choose our matrix $X(t)$
to have the following structure:
\begin{align}\label{X}
   X(t)=
     \begin{pmatrix}
       0 & x_{1,2}(t) & 0 & \cdots & 0 &0\\
       -\bar{x}_{1,2}(t) & 0 & x_{2,3}(t) & \cdots & 0 &0\\
       0 & -\bar{x}_{2,3}(t) & 0 & \cdots & 0 &0\\
       \vdots & \vdots & \ddots & \ddots & \vdots &\vdots\\
       0 & 0 & 0 & \cdots & 0 &  x_{N-1,N}(t)\\
       0 & 0 & 0 & \cdots & -\bar{x}_{N-1,N}(t) & 0 \\
     \end{pmatrix},
\end{align}
where the complex functions $x_{i,i+1}(t)$, $i=1,\cdots,N-1$, will
be chosen so that equation (\ref{eqnm1.1}) yields that $A_1$ is a
polynomial of degree $1$. From the definitions of $F$, $A_2$ and
$A_1$ (see (\ref{F}), (\ref{A2}) and (\ref{eqnm1.1})) it follows
that
\begin{align}\label{A1c}
\nonumber A_1=&((1+\alpha)I+J)J-t(J+(\alpha+2)A)+t^2(A-A^2)+\\
   &\qquad\qquad\qquad\qquad\qquad\qquad  te^{At}t^{\frac{1}{2}J}(JX(t)+X(t)J)t^{-\frac{1}{2}J}e^{-At}.
\end{align}
We now compute the right hand side of the expression above. From the
definitions of $X(t)$ and $J$ (see (\ref{X}) and (\ref{JJ})) it
follows that
\begin{align*}
   X(t)J+JX(t)=
     \begin{pmatrix}
       0 & y_{1,2}(t) & 0 & \cdots & 0 &0\\
       -\bar{y}_{1,2}(t) & 0 & y_{2,3}(t) & \cdots & 0 &0\\
       0 & -\bar{y}_{2,3}(t) & 0 & \cdots & 0 &0\\
       \vdots & \vdots & \ddots & \ddots & \vdots &\vdots\\
       0 & 0 & 0 & \cdots & 0 &  y_{N-1,N}(t)\\
       0 & 0 & 0 & \cdots & -\bar{y}_{N-1,N}(t) & 0 \\
     \end{pmatrix},
\end{align*}
where
\begin{align*}
    y_{i,i+1}(t)=(2(N-i)-1)x_{i,i+1}(t),\quad i=1,\cdots,N-1.
\end{align*}
Then, by direct calculation we have
\begin{align*}
   M(t)&=t^{\frac{1}{2}J}(X(t)J+JX(t))t^{-\frac{1}{2}J}=\\
     &\qquad =\begin{pmatrix}
       0 & t^{\frac{1}{2}}y_{1,2}(t) & \cdots & 0 &0\\
       -t^{-\frac{1}{2}}\bar{y}_{1,2}(t) & 0 &  \cdots & 0 &0\\
       0 & -t^{-\frac{1}{2}}\bar{y}_{2,3}(t)  & \cdots & 0 &0\\
       \vdots & \vdots & \ddots  & \vdots &\vdots\\
       0 & 0  & \cdots & 0 &  t^{\frac{1}{2}}y_{N-1,N}(t)\\
       0 & 0  & \cdots & -t^{-\frac{1}{2}}\bar{y}_{N-1,N}(t) & 0 \\
     \end{pmatrix}
   .
\end{align*}
Taking into account the definition of $A$ (see (\ref{AA})), this
implies that
\begin{equation*}
   \mbox{ad}_A M=
     \begin{pmatrix}
     z_{1,1}(t) & 0 & z_{1,3}(t)&\cdots & 0 &0\\
     0 & z_{2,2}(t) & 0&\cdots & 0 &0\\
       0 & 0 & z_{3,3}(t) & \cdots & 0 &0\\
       \vdots & \vdots & \vdots & \ddots & \vdots &\vdots\\
       0 & 0 & 0 & \cdots & 0 & z_{N-2,N}(t)\\
       0 & 0 & 0 & \cdots & z_{N-1,N-1}(t) &  0\\
       0 & 0 & 0 & \cdots & 0 & z_{N,N}(t) \\
     \end{pmatrix},
\end{equation*}
where (put $\nu_0=0$ and $\nu_N=0$)
\begin{equation*}
    z_{i,i}(t)=(\nu_{i-1}\bar{y}_{i-1,i}(t)-\nu_i\bar{y}_{i,i+1}(t))t^{-\frac{1}{2}},\quad
    i=1,\cdots,N,
\end{equation*}
and
\begin{equation*}
    z_{i,i+2}(t)=(\nu_{i}y_{i+1,i+2}(t)-\nu_{i+1}y_{i,i+1}(t))t^{\frac{1}{2}},\quad
    i=1,\cdots,N-2.
\end{equation*}
In order for $A_1$ to be a matrix polynomial of degree 1 we need to
take
\begin{equation*}
    y_{i,i+1}(t)=y_{i,i+1}t^{-\frac{1}{2}},
\end{equation*}
where now, abusing notation, $y_{i,i+1}\in\mathbb{C}$. We can then
write $M(t)=\frac{1}{t}Y-Y^*$, where
\begin{equation}\label{Y}
   Y=
     \begin{pmatrix}
       0 & 0 & 0 & \cdots & 0 &0\\
       -\bar{y}_{1,2} & 0 & 0 & \cdots & 0 &0\\
       0 & -\bar{y}_{2,3} & 0 & \cdots & 0 &0\\
       \vdots & \vdots & \ddots & \ddots & \vdots &\vdots\\
       0 & 0 & 0 & \cdots & 0 &  0\\
       0 & 0 & 0 & \cdots & -\bar{y}_{N-1,N} & 0 \\
     \end{pmatrix} .
\end{equation}
Hence, substituting all these expressions in (\ref{A1c}) we have
\begin{align*}
A_1 &=((1+\alpha)I+J)J+Y-t(J+(\alpha+2)A+Y^*-\mbox{ad}_A
Y)\\&\nonumber \qquad+t^2(A-A^2+\frac{1}{2}\mbox{ad}_A^2
Y-\mbox{ad}_A Y^*)+\sum _{n\ge 3}\frac{t^n}{n!}(\mbox{ad}_A^n
Y-\mbox{ad}_A^{n-1}Y^*).
\end{align*}
Since the structure of the matrices $A$ and $Y$ implies that
$\mbox{ad}_A^2 Y$ has null entries out of the diagonal $(i,i+1)$ and
that $\mbox{ad}_A Y^*$ has null entries out of the diagonal
$(i,i+2)$, $A_1$ is a matrix polynomial of degree 1 if and only if
\begin{equation}\label{Cond}
    \begin{array}{cc}
      \frac{1}{2}\mbox{ad}_A^2 Y + A&=0 \\
      \mbox{ad}_A Y^* + A^2 &=0.
    \end{array}
\end{equation}
From the first equation in (\ref{Cond}) we obtain $\displaystyle
y_{i,i+1}=-\frac{i(N-i)}{\bar{\nu}_i}$ which implies
\begin{equation}\label{xx1}
      x_{i,i+1}=-\frac{i(N-i)}{(2N-2i-1)\bar{\nu}_i},\quad
   i=1,\cdots,N-1.
\end{equation}
We are again abusing notation taking
$x_{i,i+1}(t)=x_{i,i+1}t^{-\frac{1}{2}}$, $x_{i,i+1}\in\mathbb{C}$.
From the second equation in (\ref{Cond}) we get $\displaystyle
y_{i,i+1}=\frac{\nu_i c}{\nu_1}+(i-1)\nu_i$ which implies
\begin{equation}\label{xx2}
   x_{i,i+1}=\frac{\nu_i c}{(2N-2i-1)\nu_1}+\frac{(i-1)\nu_i}{2N-2i-1},\quad
   i=1,\cdots,N-1,
\end{equation}
where $c$ is any complex number. Then, equating (\ref{xx1}) and
(\ref{xx2}) we have the following set of equations:
\begin{equation*}
    c|\nu_i|^2+(i-1)\nu_1|\nu_i|^2+i(N-i)\nu_1=0,\quad
    i=1,\cdots,N-1.
\end{equation*}
After removing the parameter $c$, we can write these equations in
the following two equivalent and more convenient ways
\begin{align}\label{prim}
    i(N-i)|\nu_{i+1}|^2&=(i+1)(N-i-1)|\nu_{i}|^2+|\nu_{i}|^2|\nu_{i+1}|^2,\\
\label{seg}
    i(N-i)|\nu_{N-1}|^2&=(N-1)|\nu_{i}|^2+(N-i-1)|\nu_{i}|^2|\nu_{N-1}|^2,
\end{align}
for $i=1,\cdots,N-2$, $N>2$.
With this choice of relations between the parameters
$\nu_1,\cdots,\nu_{N-1}$, the function $A_1$ defined in
(\ref{eqnm1.1}) is a matrix polynomial of degree 1, namely:
\begin{equation*}
    A_1(t)=((1+\alpha)I+J)J+Y-t(J+(\alpha+2)A+Y^*-\mbox{ad}_A Y).
\end{equation*}

We now show that under the assumption (\ref{prim}) (or, equivalently
(\ref{seg})), we can produce a matrix $A_0$ such that $W$ also
satisfies the differential equation (\ref{eq2.3}).

According to (\ref{eqnm1.2}), we are going to prove that there
exists a matrix $A_0$ such that the function
$$
R^{-1}(t)(F(t)A_2(t)F(t)+F'(t)A_2(t)+F(t)(t^{\alpha}e^{-t}A_2)'t^{-\alpha}e^{t}-A_0)R(t)
$$
is Hermitian. Since $R(t)=e^{At}t^{\frac{1}{2}J}U(t)$ and $U(t)$ is
unitary, it is equivalent to prove that
\begin{equation}\label{HH}
\chi(t)=t^{-\frac{1}{2}J}e^{-At}(F(t)A_2(t)F(t)+F'(t)A_2(t)+F(t)(t^{\alpha}e^{-t}A_2)'t^{-\alpha}e^{t}-A_0)e^{At}t^{\frac{1}{2}J}
\end{equation}
is always Hermitian.

We now look for a convenient expression for (\ref{HH}). We use the
formulas
\begin{align}
\nonumber    t^{-\frac{1}{2}J}e^{-At}(A_2)e^{At}t^{\frac{1}{2}J}&=t\cdot t^{-\frac{1}{2}J}Jt^{\frac{1}{2}J}=tJ,\\
\nonumber    t^{-\frac{1}{2}J}e^{-At}(A_2')e^{At}t^{\frac{1}{2}J}&=J-t^{\frac{1}{2}}A,\\
\nonumber t^{-\frac{1}{2}J}At^{\frac{1}{2}J}&=t^{-\frac{1}{2}}A,\\
\label{XX}    X(t)=t^{-\frac{1}{2}}X\quad\mbox{and}&\quad
X'(t)=-\frac{1}{2t}X(t),
\end{align}
where, again, abusing notation, $X$ is the skew-Hermitian matrix
independent of $t$ with null entries except for the diagonals
$(i,i+1)$ and $(i+1,i)$, whose entries are given in (\ref{xx1}).
From the definition of $F$ and $A_2$ (see (\ref{F}) and (\ref{A2})),
and after a straightforward computation, we obtain
\begin{align*}
\nonumber
t^{-\frac{1}{2}J}e^{-At}(FA_2F)e^{At}t^{\frac{1}{2}J}=&\frac{1}{4t}J^3+\frac{1}{2}t^{-\frac{1}{2}}(J^2A+AJ^2)
\nonumber       +AJA+tX(t)JX(t)\\
\nonumber       & +t^{\frac{1}{2}}(AJX(t)+X(t)JA)+\frac{1}{2}(J^2X(t)+X(t)J^2), \\
\label{Expr}      t^{-\frac{1}{2}J}e^{-At}(F'A_2)e^{At}t^{\frac{1}{2}J} =&-\frac{1}{2t}J^2+t^{\frac{1}{2}}(AX(t)J-X(t)AJ) \\
\nonumber       &-\frac{1}{2}X(t)J-\frac{1}{2}X(t)J^2+\frac{1}{2}JX(t)J,\\
\nonumber
t^{-\frac{1}{2}J}e^{-At}(F(t^{\alpha}e^{-t}A_2)'t^{-\alpha}e^{t})e^{At}t^{\frac{1}{2}J}
      =&\frac{\alpha+1}{2t}J^2+t^{-\frac{1}{2}}((\alpha+1)AJ-\frac{1}{2}JA)-A^2-\frac{1}{2}J^2 \\
\nonumber       & -t^{\frac{1}{2}}(AJ+X(t)A)+(\alpha+1)X(t)J-tX(t)J.
\end{align*}
Using again (\ref{XX}), we obtain an expression for the matrix
function $\chi(t)$ (see (\ref{HH})):
\begin{align*}
\chi(t)&= \frac{1}{t}\bigg[\frac{1}{4}J^2(J+2\alpha I)\bigg]+\frac{1}{2}t^{-\frac{1}{2}}\bigg[J^2A+AJ^2+J^2X+(2\alpha+1)AJ+(2\alpha+1)XJ- \\
     &\quad\quad JA+JXJ\bigg]+AJA+XJX+AJX+XJA+AXJ-XAJ-XA-\\
     &\quad\quad
     A^2-\frac{1}{2}J^2-t^{\frac{1}{2}}(AJ+XJ)-t^{-\frac{1}{2}J}e^{-At}(A_0)e^{At}t^{\frac{1}{2}J}.
\end{align*}
We define the matrix $A_0$ as follows
\begin{equation*}
    A_0=
     \begin{pmatrix}
     A_{0,1,1} &  A_{0,1,2}&0 &\cdots & 0& 0\\
     0 & A_{0,2,2} & A_{0,2,3}&\cdots & 0& 0\\
       \vdots & \vdots & \vdots & \ddots &\vdots &\vdots\\
       0 & 0 & 0 & \cdots & A_{0,N-1,N-1} &  A_{0,N-1,N}\\
       0 & 0 & 0 & \cdots & 0 & A_{0,N,N}
     \end{pmatrix}
   =V_1+V_2,
\end{equation*}
where $V_{1,i,i}=A_{0,i,i}$, $i=1,\cdots ,N$, and $V_{1,i,j}=0$
otherwise, and $V_{2}=A_{0}-V_1$. Then, we obtain:
\begin{equation*}
t^{-\frac{1}{2}J}A_0t^{\frac{1}{2}J}=V_1+t^{-\frac{1}{2}}V_2\quad\mbox{and}
\end{equation*}
\begin{equation*}
-t\cdot t^{-\frac{1}{2}J}(\mbox{ad}_A
A_0)t^{\frac{1}{2}J}=-t^{\frac{1}{2}}V_3-V_4,
\end{equation*}
where
\begin{align*}
V_{3,i,i+1}&=\begin{cases} \nu_i(A_{0,i+1,i+1}-A_{0,i,i}), & i=1,\cdots,N-1, \\
0,& \text{otherwise};\end{cases} \\
V_{4,i,i+2}&=\begin{cases} \nu_iA_{0,i+1,i+2}-\nu_{i+1}A_{0,i,i+1},& i=1,\cdots,N-2,\\
0,&\text{otherwise}.\end{cases}
\end{align*}
Taking this into account, as well as the fact that $X(JA-AJ)=XA$ and
that $\frac{1}{4}J^2(J+2\alpha I)$, $XJX$ and $-\frac{1}{2}J^2$ are
already Hermitian, in order to prove that (\ref{HH}) is Hermitian it
is enough to impose that the matrices
\begin{equation}\label{Herm1}
\frac{1}{2}[J^2A+AJ^2+J^2X+(2\alpha+1)AJ+(2\alpha+1)XJ-JA+JXJ]-V_2,
\end{equation}
\begin{equation}\label{Herm2}
-(AJ+XJ)-V_3\quad \mbox{and}
\end{equation}
\begin{equation}\label{Herm3}
(AJA+AJX+AXJ-A^2)-V_4
\end{equation}
are Hermitian and that
\begin{equation}\label{Herm4}
\mbox{ad}_A^2 A_0=0.
\end{equation}
Condition (\ref{Herm1}) allows us to define the upper diagonal
$(i,i+1)$ of the matrix $A_0$ by
\begin{equation*}
    A_{0,i,i+1}=(\alpha+N-i)[\nu_i(N-i-1)+x_{i,i+1}(2N-2i-1)],\quad
    i=1,\cdots,N-1,
\end{equation*}
which, using (\ref{xx1}) and (\ref{seg}), implies
\begin{equation}\label{A01}
A_{0,i,i+1}=-\frac{(N-1)(\alpha+N-i)\nu_i}{|\nu_{N-1}|^2},\quad
    i=1,\cdots,N-1.
\end{equation}
Condition (\ref{Herm2}) gives us a recursive expression for the the
elements $A_{0,i,i}$:
\begin{equation*}
A_{0,i,i}=A_{0,i+1,i+1}-\bigg[N-i-1+(2N-2i-1)\frac{x_{i,i+1}}{\nu_i}\bigg],\quad
i=N-1, N-2,\cdots,1.
\end{equation*}
Putting $A_{0,N,N}=0$ and using (\ref{xx1}) and (\ref{seg}) we
obtain
\begin{equation}\label{A02}
A_{0,i,i}=\frac{(N-1)(N-i)}{|\nu_{N-1}|^2},\quad
    i=1,\cdots,N.
\end{equation}
Condition (\ref{Herm3}) is then equivalent to
\begin{equation*}
    (N-i-2)\nu_i\nu_{i+1}+(2N-2i-3)\nu_ix_{i+1,i+2}+\nu_iA_{0,i+1,i+2}-\nu_{i+1}A_{0,i,i+1}=0,
\end{equation*}
which can be easily deduced using (\ref{xx1}), (\ref{A01}) and
(\ref{prim}). Finally, condition (\ref{Herm4}) is equivalent to
\begin{align*}
(\mbox{ad}_A^2
A_0)_{i,i+2}&=\nu_i\nu_{i+1}(A_{0,i,i}-2A_{0,i+1,i+1}+A_{0,i+2,i+2})=0,\\
(\mbox{ad}_A^2
A_0)_{i,i+3}&=\nu_i\nu_{i+1}A_{0,i+2,i+3}-2\nu_i\nu_{i+2}A_{0,i+1,i+2}+\nu_{i+1}\nu_{i+2}A_{0,i,i+1}=0,
\end{align*}
for $i=1,\cdots,N-2$, which can be deduced from (\ref{A01}) and
(\ref{A02}).

By looking directly at (\ref{A01}) and (\ref{A02}) it is easy to
conclude that
\begin{equation*}
    A_0=\frac{N-1}{|\nu_{N-1}|^2}[J-(\alpha I+J)A].
\end{equation*}
We have thus proved the following theorem:

\begin{thm}
Let $W$ be the weight matrix defined by (\ref{W}) where the moduli
of the entries $|\nu_i|, i=1,\cdots,N-2$, of the matrix $A$ are
defined from $\nu_{N-1}$ using the equations (\ref{seg}). Consider
the matrices $A$ and $J$ defined by (\ref{AA}) and (\ref{JJ})
respectively, and $Y$, $X(t)$ and $F(t)$ defined by
\begin{equation}\label{Y}
   Y=\begin{pmatrix}
       0 & 0 & 0 & \cdots & 0 &0\\
       \frac{N-1}{\nu_1} & 0 & 0 & \cdots & 0 &0\\
       0 & \frac{2(N-2)}{\nu_2} & 0 & \cdots & 0 &0\\
       \vdots & \vdots & \ddots & \ddots & \vdots &\vdots\\
       0 & 0 & 0 & \cdots & 0 &  0\\
       0 & 0 & 0 & \cdots & \frac{N-1}{\nu_{N-1}} & 0 \\
     \end{pmatrix},
     \end{equation}
   \begin{equation*}X(t)=t^{-\frac{1}{2}}
     \begin{pmatrix}
       0 & x_{1,2} & 0 & \cdots & 0 &0\\
       -\bar{x}_{1,2} & 0 & x_{2,3} & \cdots & 0 &0\\
       0 & -\bar{x}_{2,3} & 0 & \cdots & 0 &0\\
       \vdots & \vdots & \ddots & \ddots & \vdots &\vdots\\
       0 & 0 & 0 & \cdots & 0 &  x_{N-1,N}\\
       0 & 0 & 0 & \cdots & -\bar{x}_{N-1,N} & 0 \\
     \end{pmatrix},
\end{equation*}
\begin{equation*}
F(t)=\frac{1}{2}\bigg(\frac{1}{t}J+A\bigg)+e^{At}t^{\frac{1}{2}J}X(t)t^{-\frac{1}{2}J}e^{-At},
\end{equation*}
where $x_{i,i+1}=-\D\frac{i(N-i)}{(2N-2i-1)\bar{\nu}_i},\quad
i=1,\cdots,N-1$. Finally define the coefficients of the differential
operator $A_2$, $A_1$ and $A_0$ by:
\begin{align*}
  A_2&=t(J-At), \\
  A_1&=((1+\alpha)I+J)J+Y-t(J+(\alpha+2)A+Y^*-\mbox{ad}_A Y),\\
  A_0&=\D\frac{N-1}{|\nu_{N-1}|^2}[J-(\alpha I+J)A].
\end{align*}
Then the weight matrix  $W$ satisfies the set of equations
\begin{equation*}
\begin{array}{rcl}
A_2W&=&WA_2^* ,\\
2(A_2W)'&=&WA_1^*+A_1W,\\
(A_2W)''-(A_1W)'+A_0W&=&WA_0^*,
\end{array}
\end{equation*}
so that the second order differential operator
\begin{equation}\label{op2}
    \ell_{2,2}=D^2A_2+D^1A_1+D^0A_0
\end{equation}
is symmetric with respect to $W$. $W$ can be factorized as
$W(t)=t^{\alpha}e^{-t}R(t)R^*(t)$, $\alpha>-1$, where
$R(t)=e^{At}t^{\frac{1}{2}J}e^{2tX(t)}$ satisfies $R'(t)=F(t)R(t)$.
\end{thm}
\hfill$\Box$

We now prove that the symmetric second order operators $\ell_{2,1}$
and $\ell_{2,2}$ (see (\ref{difopeW}) and (\ref{op2})) for $W$
commute and satisfy a relation given by a polynomial of degree $N$
in two variables.

\begin{thm}
The second order differential operators $\ell_{2,1}$ and
$\ell_{2,2}$ given in (\ref{difopeW}) and (\ref{op2}), respectively,
commute, i.e. $[\ell_{2,1},\ell_{2,2}]=0$. Moreover, they satisfy
the following relation:
\begin{equation}\label{Fact}
    \prod_{i=1}^N
    \bigg((i-1)\ell_{2,1}-\ell_{2,2}+\left[\frac{(N-1)(N-i)}{|\nu_{N-1}|^2}+(i-1)(N-i)\right]I\bigg)=0.
\end{equation}
\end{thm}

\pf Taking the monic orthogonal polynomials with respect to $W$ as
eigenfunctions, there exists an isomorphism between differential
operators and their corresponding eigenvalues. Since the operators
$\ell_{2,1}$ and $\ell_{2,2}$ have a common system of
eigenfunctions, they commute if and only if their corresponding
eigenvalues commute. These eigenvalues are given by:
\begin{equation*}
    \begin{array}{rl}
      \Gamma_{n,1} & =n(A-I)+(J+\alpha I)A-J \\
      \Gamma_{n,2} & =-n^2A-n(J+(\alpha+1)A+Y^*-\mbox{ad}_A
      Y)+\D\frac{N-1}{|\nu_{N-1}|^2}(J-(\alpha I+J)A).
    \end{array}
\end{equation*}
Since $\Gamma_{n,1}$ and $\Gamma_{n,2}$ have non null entries only
in the diagonals $(i,i)$ and $(i,i+1)$, they commute if and only if
\begin{equation*}
    \begin{array}{rl}
      (\Gamma_{n,1})_{i,i+1}[(\Gamma_{n,2})_{i+1,i+1}-(\Gamma_{n,2})_{i,i}]+(\Gamma_{n,2})_{i,i+1}[(\Gamma_{n,1})_{i,i}-(\Gamma_{n,1})_{i+1,i+1}] & =0,\\
    (\Gamma_{n,1})_{i,i+1}(\Gamma_{n,2})_{i+1,i+2}-(\Gamma_{n,2})_{i,i+1}(\Gamma_{n,1})_{i+1,i+2}&=0.
    \end{array}
\end{equation*}
These formulas can be easily checked using the definitions of $A$,
$J$ and $Y$ (see (\ref{AA}), (\ref{JJ}) and (\ref{Y})) and formulas
(\ref{prim}) and (\ref{seg}). Hence, $[\ell_{2,1},\ell_{2,2}]=0$.

Now we prove (\ref{Fact}). We can use the already mentioned
correspondence between differential operators and their eigenvalues.
Let us write, for $i=1,\cdots,N$,
\begin{equation*}
    \Delta_i=(i-1)\Gamma_{n,1}-\Gamma_{n,2}+\left(\frac{(N-1)(N-i)}{|\nu_{N-1}|^2}+(i-1)(N-i)\right)I.
\end{equation*}
Then, it is enough to prove that $\Delta_1\Delta_2\cdots\Delta_N=0$.
It is easy to check that the $i$-th element of the main diagonal of
$\Delta_i$ is 0. Taking into account that $\Gamma_{n,1}$ and
$\Gamma_{n,2}$ have non null entries only in the diagonals $(i,i)$
and $(i,i+1)$, we deduce that, for $i=1,\cdots,N-1$, the $i$-th row
of $\Delta_i$ has all its entries equal to zero except the entry
$i+1$, and the row $N$ of $\Delta_N$ vanishes. Then, it is
straightforward to conclude that (\ref{Fact}) holds. \hfill $\Box$
\par\bigskip

We complete this section illustrating some structural formulas for a
sequence of orthogonal polynomials $(\mathcal P_{n,\alpha, a})_n$
with respect to the weight matrix
\begin{equation}\label{W2}
    W_{\alpha, a}(t)=t^{\alpha}e^{-t}
           \begin{pmatrix}
             t(1+|a|^2t) & at \\
             \bar{a}t & 1 \\
           \end{pmatrix}
         ,\quad a\in\mathbb{C}\setminus \{0\},\quad
         \alpha>-1,\quad t\in (0,+\infty),
\end{equation}
which is the special case of the  weight matrix $W_{\alpha
,\nu_1,\cdots,\nu_{N-1}}$ for $N=2$.

Indeed, the sequence $(\mathcal P_{n,\alpha, a})_n$ can be defined
by means of a Rodrigues' formula. Let us write
\begin{equation*}
    R_a(t)=\begin{pmatrix}
             t(1+|a|^2t) & at \\
             \bar{a}t & 1 \\
           \end{pmatrix}.
\end{equation*}
Then, the matrix polynomials defined by
\begin{align}\label{Pn}
    \nonumber\mathcal P_{0,\alpha,a}&=\begin{pmatrix}
             1&-a(1+\alpha) \\
             0 & 1 \\
           \end{pmatrix},\\
    \mathcal{P}_{n,\alpha, a}(t)&=\Phi_{n,\alpha,a}\big[t^{\alpha+n}e^{-t}(R_a(t)+X_{n,a})\big]^{(n)}R_a^{-1}(t)t^{-\alpha}e^{t},\quad
    n\geq1,
\end{align}
where
$$
\Phi_{n,\alpha,a}=\begin{pmatrix}
               1 & -a(1+\alpha) \\
               0 & 1/\gamma_{n,a} \\
             \end{pmatrix},\quad
X_{n,a}=\begin{pmatrix}
               0& -an \\
               0 & 0 \\
             \end{pmatrix},\quad \mbox{and}\quad
$$
\begin{equation}\label{gamm}
\gamma_{n,a}=1+n\vert a\vert^2,
\end{equation}
are orthogonal with respect to the weight matrix $W_{\alpha ,a}$
defined in (\ref{W2}). The leading coefficient
$\Lambda_{n,\alpha,a}$ is a nonsingular matrix given by
$$
\Lambda_{n,\alpha,a}=(-1)^n\begin{pmatrix}
                 1& -a(1+n+\alpha) \\
                0 & 1 \\
              \end{pmatrix} .
$$
This can be proved as in \cite{DG3} or \cite{DL}.

The normalization by the matrices $\Phi_{n,\alpha,a}$ allows us to
get simpler expressions for some other formulas related to
polynomials $(\mathcal P_{n,\alpha, a})_n$. For instance, one can
prove that they satisfy the following three--term recurrence
relation:
$$
t {\mathcal P}_{n,\alpha,a}(t)  = D_{n}{\mathcal
P}_{n+1,\alpha,a}(t) + E_n {\mathcal P}_{n,\alpha,a}(t) +
F_n{\mathcal P}_{n-1,\alpha,a}(t),
$$
where the matrices  $E_n$, $D_{n}$, $n\ge 0$, and $F_n$, $n\ge 1$,
are given explicitly by the expressions:
\begin{align*}
D_{n}&=-\begin{pmatrix}
 1 &a \\
 0&1
\end{pmatrix},\\
E_n&=\begin{pmatrix}
\displaystyle \frac{\gamma_{n+1,a}(2n+3+\alpha)-1}{\gamma_{n+1,a}}& \displaystyle a(1+n+\alpha) \\
\displaystyle \frac{\bar{a}}{\gamma_{n+1,a}\gamma_{n,a}} &
\displaystyle \frac{\gamma_{n,a}(2n+\alpha)+1}{\gamma_{n,a}}
\end{pmatrix},
\\
F_n&=-\D\frac{1}{\gamma_{n,a}}\begin{pmatrix}
n\gamma_{n+1,a}(1+n+\alpha)& \displaystyle 0  \\
\displaystyle \frac{n\bar{a}}{\gamma_{n,a}}
&n\gamma_{n-1,a}(n+\alpha)
\end{pmatrix}.
\end{align*}
The $L^2$-norm of $(\mathcal P_{n,\alpha, a})_n$ is given by
\begin{equation*}
\|\mathcal{P}_{n,\alpha,a}\|_{L^2(W)}^2=n!\begin{pmatrix}
                \Gamma(\alpha+n+2)\gamma_{n+1,a} & 0 \\
                0 &  \Gamma(\alpha+n+1)/\gamma_{n,a}\\
              \end{pmatrix}.
\end{equation*}
We will use the family $(\mathcal P_{n,\alpha, a})_n$ in Section
\ref{alg} to study the algebra of differential operators associated
with the weight matrix $W_{\alpha ,a}$. Before, we need some
conditions of symmetry for higher order differential operators.

\section{Conditions of symmetry for higher order differential operators}\label{Hodo}

The purpose of this section is to generalize the symmetry conditions
(\ref{eq2.1}), (\ref{eq2.2}) and (\ref{eq2.3}) of a second order
differential operator  to symmetry conditions for a higher order
differential operator with respect to a weight matrix $W$ (for the
scalar case see \cite{LK}).

Let $\ell_k$ be a differential operator of order $k$
\begin{equation}\label{dop}
\ell _k=\sum_{i=0}^k D^i A_i(t),
\end{equation}
where $D^i$ is the $i$--th derivative with respect to $t$ and
$A_i(t), i=0,\cdots,k$, are matrix polynomials of degree no bigger
than $i$. Let us write $A_i(t)$ as:
$$
A_i(t)=\sum_{j=0}^i t^jA_{j}^i ,\;\;\; A_{j}^i\in
\mathbb{C}^{N\times N}.
$$
Firstly we relate the symmetry of $\ell_k$ with the moments
$\mu_n=\int_{\Omega}t^n W(t)dt$, $n\ge 0$, of a weight matrix $W$,
where $\Omega$ is the support of $W$. In this section, $[n]_i$ will
denote the bounded factorial defined by
\begin{align}\label{BF}
[n]_i=&n(n-1)\cdots(n-i+1),\quad n\geq i>0,\\
\nonumber[n]_0=&1,\quad [n]_i=0,\quad i>n\geq0.
\end{align}

\begin{prop}
The differential operator $\ell _k$ (see (\ref{dop})) is symmetric
with respect to the weight matrix $W$ if and only if the moments
$(\mu_n)_n$ of $W$ satisfy the following $k+1$ sets of equations:
\begin{equation}\label{moment}
\sum_{i=0}^{k-l} \bigg(\begin{array}{c}
       k-i \\
       l
     \end{array}
\bigg)[n-l]_{k-l-i}B_n^{k-i}=(-1)^l(B_n^l)^*,\quad
l=0,\cdots,k,\quad n\geq l
\end{equation}
where
\begin{equation}\label{momentss}
B_n^l=\sum_{i=0}^l A_{l-i}^l\mu_{n-i},\quad l=0,\cdots,k,\quad n\geq
l.
\end{equation}
\end{prop}
\pf The symmetry of $\ell_k$ implies
\begin{equation}\label{eq1}
\int_{\Omega}\ell_k(t^n I)\,W(t)\,(t^m I)^*\,dt=\int_{\Omega} (t^n
I)\,W(t)\, \ell_k(t^m I)^*\,dt,\quad n,m\geq0.
\end{equation}
Then, from (\ref{eq1}) we have the following moment equations:
\begin{equation}\label{eq2}
\sum_{i=0}^k [n]_{k-i}B_{n+m}^{k-i}=\sum_{i=0}^k
[m]_{k-i}(B_{n+m}^{k-i})^*.
\end{equation}
In order to obtain (\ref{moment}) we use bounded complete induction
over $l$. When $l=0$ we put $m=0$ in (\ref{eq2}) and we get the
first equation.

Assume the first $j-1$ equations are true and let us show that the
$j$--th also holds. Put $n-j+1$ and $m=j-1$ in (\ref{eq2}) to get
\begin{equation}\label{eq3}
\sum_{i=0}^k [n-j+1]_{k-i}B_{n}^{k-i}=\sum_{i=0}^{j-1}
[j-1]_i(B_{n}^{i})^*.
\end{equation}
Substitute the expressions for $(B_{n}^{i})^*,\; i=0,\cdots,j-2$,
defined in (\ref{moment}), into (\ref{eq3}) and group again all the
$B_{n}^{i},\; i=0,\cdots,k$. To conclude it is enough to take into
account that for $m=0,\cdots,k$, we have that
$$
[n-j+1]_{m}=\sum_{h=0}^{m} (-1)^h \bigg(\begin{array}{c}
       m \\
       h
     \end{array}
\bigg) [j-1]_h [n-h]_{m-h}.
$$
Note that, from the definition of bounded factorial given in
(\ref{BF}), the left hand side of the formula above vanishes for
$n=j+p-1$, where $p=0,\cdots,m-1$. If we see that the right hand
side of this expression is also zero for all such values of $n$ it
will imply the equality since both expressions are polynomials in
$n$.

Write the right hand side of the formula above as
$$
n\cdots(n-m+1)\sum_{h=0}^{m} (-1)^h \bigg(\begin{array}{c}
       m \\
       h
     \end{array}
\bigg) \frac{(j-1)\cdots(j-h)}{n\cdots(n-h+1)}.
$$
Substitute $n$ by $j+p-1$, for $p=0,\cdots,m-1$ in the expression
above. Then we have to verify that
$$
(j+p-1)\cdots(j+p-m)\sum_{h=0}^{m} (-1)^h \bigg(\begin{array}{c}
       m \\
       h
     \end{array}
\bigg) \frac{(j-1)\cdots(j-h)}{(j+p-1)\cdots(j+p-h)}=0.
$$
But this is always true considering the fact that $\sum_{h=0}^{m}
(-1)^h \bigg(\begin{array}{c}
       m \\
       h
     \end{array}
\bigg) h^p=0$ for $p=0,\cdots,m-1$ and that the formula above is
always a polynomial in $h$ of degree at most $m-1$.

The converse is similar.\hfill $\Box$

Using the previous proposition, we can guarantee the symmetry of $\ell
_k$ with respect to $W$ from a set of differential equations which
relates $W$ with the  coefficients of the differential operator $\ell _k$.

\begin{thm}\label{scod}
Assume the boundary conditions that
\begin{equation}\label{bound}
\sum_{i=0}^{p-1}(-1)^{k-i+p-1} \bigg(\begin{array}{c}
       k-i \\
       l
     \end{array}
\bigg) \big(A_{k-i}\cdot W\big)^{(p-1-i)},\quad p=1,\cdots,k,\quad
l=0,\cdots,k-p
\end{equation}
should have vanishing limits at each of the endpoints of the support
of $W$ and that the following $k+1$
equations are satisfied:
\begin{equation}\label{symeqs}
\sum_{i=0}^{k-l}(-1)^{k-i-l} \bigg(\begin{array}{c}
       k-i \\
       l
     \end{array}
\bigg) \big(A_{k-i}\cdot W\big)^{(k-i-l)}=(-1)^l W\cdot A_l^*,\quad
l=0,\cdots,k.
\end{equation}
Then the differential operator $\ell_k$ (see (\ref{dop})) is
symmetric with respect to $W$.

\end{thm}
\pf Using the moment equations (\ref{momentss}) we obtain

\begin{eqnarray*}
\lefteqn{[n-l]_{k-l-i}B_n^{k-i}=[n-l]_{k-l-i}\bigg(\sum_{j=0}^{k-i}A_{k-i-j}^{k-i}\mu_{n-j}\bigg)=}\\
\noalign{\smallskip} &
&\;\;\quad\qquad\qquad=\int_{\Omega}[n-l]_{k-l-i}t^{n-k+i}A_{k-i}(t)\cdot
W(t) dt.
\end{eqnarray*}
Integrating by parts $k-i-l$ times ($i\leq k-l-1$) the previous
integral is equal to
\begin{align*}
&\sum_{j=1}^{k-i-l}(-1)^{j-1}[n-l]_{k-l-i-j}t^{n-k+i+j}\big(A_{k-i}\cdot
W\big)^{(j-1)}\bigg]_{\partial\Omega}+\\
&\qquad\qquad\qquad\qquad \qquad\qquad\qquad\qquad
(-1)^{k-i-l}\int_{\Omega}t^{n-l}\big(A_{k-i}\cdot
W\big)^{(k-i-l)}dt.
\end{align*}
Replacing the value of $[n-l]_{k-l-i}B_n^{k-i}$ in (\ref{moment}), we
have to prove that
\begin{align*}
\sum_{i=0}^{k-l-1}&\bigg(\begin{array}{c}
       k-i \\
       l
     \end{array}
\bigg)\bigg(\D\sum_{j=1}^{k-i-l}(-1)^{j-1}[n-l]_{k-l-i-j}t^{n-k+i+j}\big(A_{k-i}\cdot
W\big)^{(j-1)}\bigg]_{\partial\Omega}\bigg)\\
&+\int_{\Omega}t^{n-l}\bigg(\sum_{i=0}^{k-l}(-1)^{k-i-l}
\bigg(\begin{array}{c}
       k-i \\
       l
     \end{array}
\bigg) \big(A_{k-i}\cdot W\big)^{(k-i-l)}-(-1)^l W\cdot
A_l^*\bigg)dt=0.
\end{align*}
But
\begin{enumerate}
  \item $\D\sum_{i=0}^{k-l-1}\bigg(\begin{array}{c}
       k-i \\
       l
     \end{array}
\bigg)\bigg(\D\sum_{j=1}^{k-i-l}(-1)^{j-1}[n-l]_{k-l-i-j}t^{n-k+i+j}\big(A_{k-i}\cdot
W\big)^{(j-1)}\bigg]_{\partial\Omega}\bigg)=0$,\\ for all
$n-l\geq0$. This follows expanding the sum, grouping in terms of the
powers of $t$, taking into account the boundary conditions
(\ref{bound}) and the fact that $W$ must have finite moments.
  \item $\D\int_{\Omega}t^{n-l}\bigg(\sum_{i=0}^{k-l}(-1)^{k-i-l}
\bigg(\begin{array}{c}
       k-i \\
       l
     \end{array}
\bigg) \big(A_{k-i}\cdot W\big)^{(k-i-l)}-(-1)^l W\cdot
A_l^*\bigg)dt=0,$\\for all $n-l\geq0$, because of the symmetry
equations (\ref{symeqs}).
\end{enumerate}
\hfill $\Box$

\bigskip

For $k=2$, Theorem \ref{scod} gives the conditions (\ref{eq2.1}),
(\ref{eq2.2}) and (\ref{eq2.3}) for the symmetry of a second order
differential operator.

For $k=3$ those conditions are
\begin{align}\nonumber
A_3W+WA_3^*&=0 ,\\\nonumber -3(A_3W)'+A_2W&=WA_2^*,\\\label{seo3}
-3(A_3W)''+2(A_2W)'&=WA_1^*+A_1W,\\ \nonumber
-(A_3W)'''+(A_2W)''-(A_1W)'+A_0W&=WA_0^*,
\end{align}
and
\begin{align*}
    A_3W,\quad 3(A_3W)'-2(A_2W),\quad -(A_3W)'+(A_2W),\quad(A_3W)''-(A_2W)'+A_1W
\end{align*}
should vanish at the endpoints of the support of $W$.

\section{The algebra of differential operators associated with $W_{\alpha,a}$}\label{alg}

The main goal of this section is to study the structure of the
algebra $\mathcal{D}(W_{\alpha ,a})$ of differential operators
having the family $(\mathcal P_{n,\alpha ,a})_n$ of orthogonal
polynomials with respect to $W_{\alpha ,a}$ (given in (\ref{W2})) as
common eigenfunctions:
$$
\mathcal{D}(W_{\alpha ,a})=\left\{\ell =\sum_{i=0}^k D^i A_i(t) :
\ell (\mathcal
P_{n,\alpha,a}(t))=\Gamma_n(\ell)\mathcal{P}_{n,\alpha
,a}(t),\;n\geq0\right\},
$$
where $\Gamma_n (\ell )$ does not depend on $t$.

The product in $\mathcal{D}(W_{\alpha ,a})$ is the composition of
operators in the usual way, which it is not commutative. If $\ell_r,
\ell_s\in \mathcal{D}(W_{\alpha ,a})$ are of order $r$ and $s$
respectively, the operator $\ell_r\ell_s$ is of order less than or
equal to $r+s$, but not necessarily $r+s$, since the leading
coefficients of both operators can be singular.
$\mathcal{D}(W_{\alpha ,a})$ is both a complex vector space and an
algebra under composition which is associative and distributive with
respect to the field we are working.

The operators in $\mathcal{D}(W_{\alpha ,a})$ need not to be
symmetric with respect to $W_{\alpha ,a}$. However, for every non
symmetric operator $\ell \in \mathcal{D}(W_{\alpha ,a})$, we can
associate an adjoint operator $\ell ^* \in \mathcal{D}(W_{\alpha
,a})$ so that $\ell +\ell ^*$ is symmetric with respect to
$W_{\alpha ,a}$. The previous assertion is a general fact for the
algebra of operators associated with any weight matrix $W$. As a
consequence, each operator $\ell \in \mathcal{D}(W_{\alpha ,a})$ can
be written uniquely in the form $\ell =\ell _1+\imath\ell_2$ for
certain symmetric operators $\ell_1, \ell _2$ (see \cite{GT} for
details).

Using the fact that the leading coefficient of $\mathcal P_n(t)$ (to
simplify the notation we remove the dependence on $\alpha $ and $a$)
is a non singular matrix, we can see that each
coefficient $A_i(t)$, $i=1, \cdots , k$, of a differential operator $\ell$ in
$\mathcal{D}(W)$ has to be a  matrix polynomial of degree less than
or equal to $i$. We observe that the map between differential
operators and the corresponding eigenvalues given by
$$
\Gamma_n:\mathcal{D}(W)\longrightarrow \mathbb{C}^{2\times 2},\quad
n=0,1,2,\ldots
$$
is a faithful representation. The property
$\Gamma_n(\ell_1\ell_2)=\Gamma_n(\ell_1)\Gamma_n(\ell_2)$ with
$\ell_1,\ell_2\in\mathcal{D}(W)$ is easy to show as well as the fact
that if $\Gamma_n(\ell)=0$, $n\geq0$, then $\ell =0$.

It is important to note that the algebra $\mathcal{D}(W)$ is
independent from the choice of the family of  orthogonal matrix
polynomials. The eigenvalues are changed by an $n$ dependent
conjugation. The family we use in this section is the sequence of
orthogonal polynomials $(\mathcal P_n)_n$ defined by the Rodrigues'
formula (\ref{Pn}).

In the following table, obtained by direct computations, we exhibit
the number of  new linearly independent differential operators
(modulo operators of lower order) that appear as one increases the
order of the operators in question:

\begin{center}\begin{tabular}{|c||c|c|c|c|c|c|c|c|c|c|}
  \hline order & 0 & 1 & 2 & 3 & 4 & 5 & 6 & 7 & 8\\\hline
  dimension & 1 & 0 & 2 & 2 & 2 & 2 & 2 & 2 & 2\\ \hline
\end{tabular}\end{center}

There are no first order differential operators. The system
$\ell_{2,1}$ and $\ell_{2,2}$ of two linearly independent
(symmetric) second order differential operators found in Section
\ref{sect2}, namely
\begin{align*}
\ell_{2,1}=D^2tI+&D^1\begin{pmatrix}
                     \alpha+2-t & at \\
                      0 & \alpha+1-t \\
                     \end{pmatrix}
                     +D^0\begin{pmatrix}
                         -1 & (1+\alpha)a \\
                            0 & 0 \\
                           \end{pmatrix},\\
\ell_{2,2}=D^2\begin{pmatrix}
                     t & -at^2 \\
                     0 & 0 \\
                   \end{pmatrix}
    +&D^1\begin{pmatrix}
         \alpha+2 & -\D\frac{1}{\bar{a}}(1+|a|^2(\alpha+2))t \\
          \D\frac{1}{a} & -t \\
          \end{pmatrix}
         +D^0\begin{pmatrix}
             \D\frac{1}{|a|^2} & -\D\frac{1+\alpha}{\bar{a}} \\
             0 & 0 \\
             \end{pmatrix},
\end{align*}
is a basis for the operators in $\mathcal{D}(W_{\alpha ,a})$ of
order $2$. Their respective eigenvalues are given by
$$
\Gamma _{n,2,1}=\begin{pmatrix}-n-1 & 0 \\ 0 & -n
\end{pmatrix},\quad \Gamma _{n,2,2}=\begin{pmatrix}\D\frac{1}{|a|^2}
& 0\\ 0 & -n \end{pmatrix}.
$$

There are two linearly independent (symmetric) third order
differential operators (modulo operators of lower order). As we
explained in the introduction, this is a phenomenon which does not
happen in the scalar case. In fact, this is the first example of
weight matrix which does not reduce to scalar weights having
symmetric differential operators of odd order.

A basis for the operators of order $3$ in $\mathcal{D}(W_{\alpha
,a})$ (modulo operators of lower order) is given by the operator
$\ell_{3,1}$ defined in (\ref{ope3}) and by $\ell_{3,2}$ defined as
\begin{align*}
\ell_{3,2}=&D^3\begin{pmatrix}
   |a|^2t^2 & at^2(-1+|a|^2t) \\
   \bar{a}t & -|a|^2t^2 \\
   \end{pmatrix}+\\
   &D^2\begin{pmatrix}
    |a|^2t(\alpha+5) & -at(-2\alpha-4+t(3+|a|^2(\alpha+5))) \\
    \bar{a}(\alpha+2) & -|a|^2t(\alpha+2) \\
    \end{pmatrix}+\\
    &D^1\begin{pmatrix}
     2|a|^2(\alpha+2)+t & a(\alpha+1)(\alpha+2)-t\bigg(\D\frac{1}{\bar{a}}+2a(2+|a|^2)(\alpha+2)\bigg) \\
     -\D\frac{1}{a} & -t \\
      \end{pmatrix}+
       \\
&D^0\begin{pmatrix}
     1+\alpha & -\D\frac{1}{\bar{a}}(1+\alpha)(1+|a|^2(\alpha+2)) \\
    \D\frac{1}{a} & -(1+\alpha) \\
     \end{pmatrix}.
\end{align*}
Their respective eigenvalues are given by (recall the definition of
$\gamma_{n,a}$ given in (\ref{gamm})):
$$
\Gamma _{n,3,1}=\frac{1}{|a|^2}\begin{pmatrix}0& \displaystyle
a(1+\alpha+n)\gamma_{n,a}\gamma_{n+1,a}\\
\bar{a}&0 \end{pmatrix},\; \Gamma
_{n,3,2}=\frac{1}{|a|^2}\begin{pmatrix}0& \displaystyle
-a(1+\alpha+n)\gamma_{n,a}\gamma_{n+1,a}\\
\bar{a}&0 \end{pmatrix}.
$$
The operator $\ell _{3,1}$ is symmetric and hence satisfies the
third order symmetry equations given by (\ref{seo3}) with their
corresponding boundary conditions. The operator $\ell_{3,2}$ is not
symmetric but skew--symmetric so that $\imath\ell_{3,2}$ is
symmetric.

The leading coefficients $A_{3,i}(t)$, $i=1,2$, of both operators
can be written as
\begin{align}\label{A3}
A_{3,1}(t)=t\bigg(At-\sum_{n\geq0}\frac{t^n}{n!}\mbox{ad}_A^n
    A^*\bigg),\quad A_{3,2}(t)=t\bigg(At+\sum_{n\geq0}\frac{t^n}{n!}\mbox{ad}_A^n
    A^*\bigg).
\end{align}
From those expressions it is easy to verify (using
$t^{-\frac{1}{2}J}At^{\frac{1}{2}J}=t^{-\frac{1}{2}}A$,
$t^{-\frac{1}{2}J}A^*t^{\frac{1}{2}J}=t^{\frac{1}{2}}A^*$ and
$\sum_{n\geq0}\frac{t^n}{n!}\mbox{ad}^n_AA^*e^{At}=e^{At}A^*$) that
\begin{align*}
   t^{-\frac{1}{2}J}e^{-At}(A_{3,1})e^{At}t^{\frac{1}{2}J}=t^{-\frac{1}{2}J}(At^2-tA^*)t^{\frac{1}{2}J}=t^{\frac{3}{2}}(A-A^*),
\end{align*}
implying that it verifies the first of the third order symmetry
equations (\ref{seo3}) since $t^{\frac{3}{2}}(A-A^*)$ is skew-
Hermitian. An analogous argument is valid for $A_{3,2}$.

In order to propose a candidate for a system of operators which
generate the full algebra $\mathcal D(W)$ and have easier
expressions for relations between generators, let us introduce a
different basis for the differential operators of order $2$:
\begin{align*}
L_1&=\ell_{2,1}-\D\frac{1}{|a|^2}I,\\
L_2&=2\ell_{2,2}-\ell_{2,1}-\D\frac{1}{|a|^2}I,
\end{align*}
while for the operators of order 3 we write $L_3=\ell_{3,1}$ and
$L_4=\ell_{3,2}$. The corresponding system of eigenvalues associated
with each operator is given by (recall again the definition of
$\gamma_{n,a}$ given in (\ref{gamm}))
\begin{align*}
\Gamma _{n,1}&=-\frac{1}{|a|^2}\begin{pmatrix} \gamma_{n+1,a} & 0 \\ 0 & \gamma_{n,a}\end{pmatrix},\\
\Gamma _{n,2}&=\frac{1}{|a|^2}\begin{pmatrix} \gamma_{n+1,a} & 0 \\ 0 & -\gamma_{n,a}\end{pmatrix},\\
\Gamma _{n,3}&=\frac{1}{|a|^2}\begin{pmatrix}0& \displaystyle
a(1+\alpha+n)\gamma_{n,a}\gamma_{n+1,a}\\
\bar{a}&0 \end{pmatrix},\\
\Gamma _{n,4}&=\frac{1}{|a|^2}\begin{pmatrix}0& \displaystyle
-a(1+\alpha+n)\gamma_{n,a}\gamma_{n+1,a}\\
\bar{a}&0 \end{pmatrix}.
\end{align*}
The set $\{I,L_1,L_2,L_3,L_4\}$ is linearly independent and from the
previous expressions for their eigenvalues it follows easily that
they satisfy a number of relations, namely:
\begin{enumerate}
\item Four quadratic relations:
\begin{align}\label{rel1}
 \nonumber L_1^2&=L_2^2, \\
  L_3^2&=-L_4^2, \\
 \nonumber L_1L_2&=L_2L_1, \\
 \nonumber L_3L_4&=-L_4L_3.
\end{align}
\item Four permutational relations:
\begin{align}\label{rel2}
 \nonumber L_1L_3-L_2L_4&=0, \\
  L_2L_3-L_1L_4&=0, \\
 \nonumber L_3L_2+L_4L_1&=0, \\
 \nonumber L_3L_1+L_4L_2&=0.
\end{align}
\item Four more second degree relations:
\begin{align}\label{rel3}
 \nonumber L_3&=L_1L_4-L_4L_1, \\
  L_4&=L_1L_3-L_3L_1,\\
 \nonumber L_3&=L_2L_3+L_3L_2, \\
 \nonumber L_4&=L_2L_4+L_4L_2.
\end{align}
\item Finally, we can present two interesting third degree relations:
\begin{align}\label{rel4}
  L_1L_3^2&=L_3^2L_1, \\
\nonumber L_2L_3^2&=L_3^2L_2.
\end{align}
\end{enumerate}
The second equation in (\ref{rel3}) allows us to remove the operator
$L_4$ as a generator of the algebra $\mathcal D(W)$.

Computational evidences allows us to conclude that a possible basis
for the operators of even order $2k$ in $\mathcal D(W)$ (modulo
operators of lower order) is given by $L_1^k$ and $L_1^{k-1}L_2$,
while for operators of odd order one can take $L_1^kL_3-L_3L_1^k$
and $L_2^kL_3+L_3L_2^k$ for $4k-1, k\geq1$, and $L_1^kL_3+L_3L_1^k$
and $L_2^kL_3-L_3L_2^k$ for $4k+1, k\geq1$.

We finally present one more relation which allows us to conjecture
that the full algebra of differential operators is generated by
$\{I, L_1, L_3\}$. Indeed,
\begin{align*}
\left[ |a|^2(2+\alpha)-1\right]&\left[|a|^2(\alpha-1)-1\right]L_2=
2|a|^2\left[ |a|^2(2\alpha +1)-2\right]L_1\\&\quad+\left[
|a|^4(\alpha^2+\alpha -5)-|a|^2(2\alpha+1)+1\right]L_1^2\\&\quad
-2|a|^2\left[ |a|^2(2\alpha +1)-2\right]L_1^3+3|a|^4L_1^4\\
&\quad -\frac{1}{2}\left[ |a|^2(2\alpha
+1)-2\right]L_3^2+\frac{15}{2}|a|^2L_3^2L_1-\frac{9}{2}|a|^2L_3L_1L_3.
\end{align*}
Note that for the exceptional values of $\alpha=1+\D\frac{1}{|a|^2}$
or $\alpha=-2+\D\frac{1}{|a|^2}$, the left side of the previous
formula vanishes. For these cases, computational evidences allow us
to conjecture that the algebra is generated by $\{I, L_1, L_2,
L_3\}$.

If one is tempted to study the algebra of differential operators for
higher matrix dimensions, the following table illustrates the number
of linearly independent differential operators that appear as one
increases the matrix dimension and order of the  differential
operators:

\begin{center}\begin{tabular}{c|c|c|c|c|c|c|c|c|c|c|}
  & $k=0$ & $k=1$ & $k=2$ & $k=3$ & $k=4$ & $k=5$ & $k=6$ & $k=7$ & $k=8$\\\hline
  $N=2$ & 1 & 0 & 2 & 2 & 2 & 2 & 2 & 2 & 2\\\hline
   $N=3$ & 1 & 0 & 2 & 0 & 3 & 4 & 5 & 6 & 16\\\hline
   $N=4$ & 1 & 0 & 2 & 0 & 3 & 0 & 4 & 6 & 9\\
\end{tabular}\end{center}

It is interesting that for $N=3, 4$, there are no third order
differential operators, but there are fifth and seventh order,
respectively. The reason could be that the coefficients (\ref{A3})
for these cases are matrix polynomials of degree at least 4 (due to
the fact that $A$ is a nilpotent matrix of order $N$).




\begin{thebibliography}{00}




\bibitem{CMV} Cantero, M. J., Moral, L. and Velázquez, L.,
\textit{Matrix orthogonal polynomials whose derivatives are also
orthogonal}, J. Approx. Theory \textbf{146} (2007), 174--211.

\bibitem{CG1}\textrm{Castro, M. M. and Grünbaum, F. A.},
\textit{Orthogonal matrix polynomials satisfying first order
differential equations: a collection of instructive examples}, J.
Nonlinear Math. Physics $\mathbf{12}$ (2005), supplement 2, 63--76.

\bibitem{CG2} Castro, M. M. and Gr\"unbaum, F. A.,
\textit{The algebra of differential operators associated to a given
family of matrix valued orthogonal polynomials: five instructive
examples}, International Math. Research Notices $\mathbf{2006}$
(2006), Article ID 47602, 33 pages.


\bibitem{D2} Durán, A. J.,
\textit{How to find weight matrices having symmetric second order
differential operators with matrix leading coefficient}, submitted.


\bibitem{DG1} Durán, A. J. and  Gr\"unbaum, F. A., {\it Orthogonal
matrix polynomials satisfying second order differential equations},
International Math. Research Notices, {\bf 2004 : 10} (2004),
461--484.


\bibitem{DG3}  Durán, A. J.  and Gr\"unbaum, F. A., {\it Structural formulas for orthogonal
matrix polynomials satisfying second-order differential equations,
I}, Constr. Approx. {\bf 22} (2005), 255--271.

\bibitem{DG2} Durán, A. J. and  Gr\"unbaum, F. A., {\it Matrix differential equations and
scalar polynomials satisfying higher order recursions}, submitted.

\bibitem{DL} Durán, A. J. and L{\'o}pez-Rodr{\'\i}guez, P.,   {\it
Structural formulas for orthogonal matrix polynomials satisfying
second order differential equations, II}, Constr. Approx.
\textbf{26} (2007), 29--47.

\bibitem{dIG} Gr\"unbaum, F. A. and de la Iglesia, M. D.,
\textit{Matrix valued orthogonal polynomials related to SU$(N+1)$,
their algebras of differential operators and the corresponding
curves}, Exp. Math. $\mathbf{16}$, No. 2, (2007), 189--207.


\bibitem{GPT1} Gr\"unbaum, F. A., Pacharoni, I., and Tirao, J. A., {\em
Matrix valued spherical functions associated to the complex
projective plane}, J. Funct. Anal. {\bf 188} (2002) 350--441.

\bibitem{GPT2} Gr\"unbaum, F. A., Pacharoni, I., and Tirao, J.
A., {\em Matrix valued orthogonal polynomials of the Jacobi type},
Indag. Mathem. {\bf 14} 3,4 (2003), 353--366.


\bibitem{GPT3} Gr\"unbaum F. A., Pacharoni I. and Tirao J.
A. {\em Matrix valued orthogonal polynomials of the Jacobi type: The
role of group representation theory}, Ann. Inst. Fourier, Grenoble,
{\bf 55} , 6, (2005), 2051--2068.


\bibitem{GT} Gr\"unbaum, F. A. and Tirao, J. A., {\em
The algebra of differential operators associate to a weight matrix},
Integr. Equ. Oper. Theory (2007), to appear.

\bibitem{K}\textrm{Krall, H. L.},
\textit{Certain differential equations for Tchebycheff polynomials},
Duke. Math. $\mathbf{4}$ (1938), 705--718.

\bibitem{L}\textrm{Littlejohn, L. L.},
\textit{Symmetry factors for differential equations},
Amer. Math. Monthly $\mathbf{7}$ (1983), 462--464.

\bibitem{LK} Littlejohn, L.L. and Krall, A.M.
\textit{Orthogonal polynomials and higher order singular Sturm-Liouville systems},
Acta Appl. Math. {\bf 17}  (1989), 99--170.

\bibitem{M1} Miranian, L., {\em On classical orthogonal polynomials and
differential operators}, J. Phys. A: Math. Gen. {\bf 38} (2005)
6379--6383.

\end{thebibliography}
\end{document}